\documentclass{amsart}
\usepackage{graphicx}
\usepackage{amsfonts, amsmath, amsfonts, amssymb}
\newtheorem{thm}{Theorem}[section]

\newtheorem{lem}[thm]{Lemma}
\newtheorem{prop}[thm]{Proposition}
\theoremstyle{definition}
\newtheorem{defn}[thm]{Definition}
\theoremstyle{remark}
\newtheorem{rem}[thm]{Remark}
\numberwithin{equation}{section}
\newenvironment{prf}{\noindent{\bf Proof}}{\\ \hspace*{\fill}$\Box$ \par}

\begin{document}

\title[The Existence of Quasimeromorphic Mappings]{The Existence of Quasimeromorphic Mappings}%
\author{Emil Saucan}%
\address{Department of Mathematics, Technion, Haifa and Software Engineering Department, Ort Braude College, Karmiel, Israel}%
\email{semil@tx.technion.ac.il}%

\thanks{This paper represents part of the authors Ph.D. Thesis written under the supervision of Prof. Uri Srebro}%

\date{27.4.2003.}
\begin{abstract}
We prove that a Kleinian group $G$ acting upon $\mathbb{H}^{n}$ admits a non-constant $G$-automorphic function,
even if it has torsion elements, provided that the orders of the elliptic (torsion) elements are uniformly bounded. 
This is accomplished by developing a technique for mashing distinct fat triangulations while preserving fatness.
\end{abstract}
\maketitle

\section{Introduction}
The object of this article 
is the study of the existence of $G$-automorphic quasimeromorphic mappings $f:\mathbb{H}^n \rightarrow
\widehat{\mathbb{R}^n}$,  $\widehat{\mathbb{R}^n} = \mathbb{R}^n \bigcup \,\{\infty\}$; i.e. such that
 \begin{equation} f(g(x))= f(x) \quad ;
  \qquad \forall x \in \mathbb{H}^n ; \forall g \in G \; ;
 \end{equation}
were $G$ is Kleinian group acting upon $\mathbb{H}^n$, and where quasimeromorphic mappings are defined as follows:
\begin{defn}
Let $D \subseteq \mathbb{R}^n$ be a domain; $n \geq 2$, and let $f: D \rightarrow \mathbb{R}^m$.
\\ $f$ is called $ACL$ ({\it absolutely continuous on lines}) iff:
\\ (i) $f$ is continuous
\\ (ii) for any $n$-interval $Q = \overline{Q} = \{a_i \leq x_i \leq b_i\,|\, i=1,\ldots,n\}$, $f$ is absolutely
continuous on almost every line segment in $Q$, parallel to the coordinate axes.
\end{defn}

\begin{lem} [\cite{v}, 26.4]
If $f:D \subseteq \mathbb{R}^n \rightarrow  \mathbb{R}^m$ is $ACL$, then $f$ admits partial derivatives almost
everywhere.
\end{lem}

The result above justifies the following Definition:

\begin{defn}
$f:D \subseteq \mathbb{R}^n  \rightarrow  \mathbb{R}^m$ is $ACL^p$ iff its derivatives are locally $L^p$
integrable, $p \geq 1$.
\end{defn}

\begin{defn}
Let $D \subseteq \mathbb{R}^n$ be a domain; $n \geq 2$ and let $f: D \rightarrow \mathbb{R}^m$ be a continuous
mapping. $f$ is called
\begin{enumerate}
\item {\it quasiregular} iff (i) $f$ is $ACL^n$ and
\\ \hspace*{2.3cm}(ii) $\exists K \geq 1$ s.t.
\begin{equation}
|f'(x)| \leq KJ_f(x)\; a.e.
\end{equation}
\\ where $f'(x)$ denotes the formal derivative of $f$ at $x$, $|f'(x)| = \sup \raisebox{-0.2cm}{\mbox{\hspace{-0.7cm}\tiny$|h|=1$}}|f'(x)h|$, and where $J_f(x) = detf'(x)$.
\\ The smallest $K$ that satisfies (4.1) is called the {\it outer dilatation} of $f$.
\item {\it quasiconformal} iff $f:D \rightarrow f(D)$ is a local homeomorphism.
\item {\it quasimeromorphic} iff $f:D \rightarrow \widehat{\mathbb{R}^n}$, $\widehat{\mathbb{R}^n} = \mathbb{R}^n \bigcup
\,\{\infty\}$ is quasiregular, where the condition of quasiregularity at $f^{-1}(\infty)$ can be checked by
conjugation with auxiliary M\"{o}bius transformations.
\end{enumerate}
\end{defn}

\begin{rem}
One can extend the definitions above to oriented $\mathcal{C}^\infty$ Riemannian $n$-manifolds by using coordinate
charts.
\end{rem}

Our principal goal is to prove:

\begin{thm} Let $G$ be a Kleinian group with torsion acting upon $\mathbb{H}^n, \, n \geq 3$. \\ If the elliptic elements (i.e. torsion elements) of $G$ have
uniformly bounded orders,
\\then there exists a non constant $G$-automorphic quasimeromorphic mapping
\\$f: \mathbb{H}^n \rightarrow \widehat{\mathbb{R}^n}.$
\end{thm}

The question whether quasimeromorphic mappings exist was originally posed by Martio and Srebro in \cite{ms1}\,;
subsequently in \cite{ms2} they proved the existence of fore-mentioned mappings in the case of co-finite groups
i.e. groups such that $ Vol_{hyp} (H^n/ G) < \infty  $ (the important case of geometrically finite groups being
thus included). Also, it was later proved by Tukia  (\cite{tu}) that the existence of non-constant
quasimeromorphic mappings (or qm-maps, in short) is assured in the case when G acts torsionless upon
$\mathbb{H}^n$. Moreover, since for torsionless Kleinian groups $G$, $\mathbb{H}^n/G$ is a (analytic) manifold,
the next natural question to ask is whether there exist qm-maps $f:M^n \rightarrow \widehat {\mathbb{R}^n} $;
where $M^n$ is an orientable $n-$manifold. The affirmative answer
to this question is due to K.
 Peltonen (see \cite{pe}); to
be more precise she proved the existence of qm-maps in the case when $M^n$ is a connected, orientable $
\mathcal{C^\infty}$\!-Riemannian manifold.
 \newline In contrast with the above results it was proved by Srebro (\cite{sr}) that, for any $n\geq 3$, there exists a Kleinian group $G \;\rhd \! \!\! \!<
 \mathbb{H}^n$ s.t. there exists no non-constant, $G-$automorphic function  $f:\mathbb{H}^n\rightarrow \mathbb{R}^n$.
 More precisely, if $G$ (as above) contains elliptics of unbounded orders (with non-degenerate fixed set), then $G$ admits no non-constant
$G-$automorphic qm-mappings.
\newline Since all the existence results were obtained in constructive manner by using the classical "Alexander
trick" (see \cite{al}), it is only natural that we try to attack the problem using the same method. For this
reason we present here in succinct manner Alexander's method: One starts by constructing a suitable triangulation
(Euclidian or hyperbolic) of $\mathbb{H}^n$. Since $\mathbb{H}^n$ is orientable, an orientation consistent with
the given triangulation  (i.e. such that two given $n$-simplices having a $(n-1)$-dimensional face in common will
have opposite orientations) can be chosen. Then one quasiconformally maps the simplices of the triangulation into
$\widehat{\mathbb{R}^n}$ in a chess-table manner: the positively oriented ones onto the interior of the standard
simplex in $\mathbb{R}^n$ and the negatively oriented ones onto its exterior.  To ensure the existence of such a
chessboard triangulation, a further barycentric type of subdivision may be required, rendering a triangulation
whose simplices satisfy the condition that every $(n-2)$-face is incident to an even number of $n$-simplices. If
the dilatations of the $qc$-maps constructed above are uniformly bounded, then the resulting map will be
quasimeromorphic.
\newline The dilatations of each of the $qc$-maps above is dictated by the proportions of the respective simplex (see \cite{tu} , \cite{ms2}), and since the dilatation is to
be uniformly bounded, we are naturally  directing our efforts in the construction of a {\bf fat} triangulation,
where:

\begin{defn} Let $\tau \subset \mathbb{R}^n$ ; $0 \leq k \leq n$ be a $k$-dimensional simplex.
The {\it fatness}  $\varphi$ of $\tau$ is defined as being:
\begin{equation}
\varphi = \varphi(\tau) = \inf_{\hspace{0.4cm}\sigma < \tau
\raisebox{-0.25cm}{\hspace{-0.9cm}\mbox{\scriptsize$dim\,\sigma=l$}}}\!\!\frac{Vol(\sigma)}{diam^{l}\,\sigma}
\end{equation}
The infimum is taken over all the faces of $\tau$, $\sigma < \tau$, and $Vol_{eucl}(\sigma)$ and $diam\,\sigma$
stand for the Euclidian $l$-volume and the diameter of $\sigma$ respectively. (If $dim\,\sigma = 0$, then
$Vol_{eucl}(\sigma) = 1$, by convention.)
\\ A simplex $\tau$ is $\varphi_0${\it -fat}, for some $\varphi_0 > 0$, if $\varphi(\tau) \geq \varphi_0$. A triangulation (of a submanifold of $\mathbb{R}^n$) $\mathcal{T} = \{ \sigma_i \}_{i\in \bf I}$ is
$\varphi_0${\it -fat} if all its simplices are $\varphi_0$-fat. A triangulation $\mathcal{T} = \{ \sigma_i
\}_{i\in \bf I }$ is {\it fat} if there exists $\varphi_0 > 0$ such that all its simplices are
$\varphi_0${\it-fat}.
\end{defn}

\begin{rem} There exists a constant $c(k)$ that depends solely upon the dimension $k$ of  $\tau$ s.t.
\begin{equation}
\frac{1}{c(k)}\cdot\varphi(\tau) \leq \min_{\hspace{0.4cm}\sigma < \tau
\raisebox{-0.25cm}{\hspace{-0.9cm}\mbox{\scriptsize$dim\,\sigma=l$}}}\hspace{-0.3cm}\measuredangle(\tau,\sigma)
\leq c(k)\cdot\varphi(\tau)\,,
\end{equation}
and
\begin{equation}
\varphi(\tau) \leq \frac{Vol(\sigma)}{diam^{l}\,\sigma} \leq c(k)\cdot\varphi(\tau)\,;
\end{equation}
where $\measuredangle(\tau,\sigma)$ denotes the  ({\it internal}) {\it dihedral angle} of $\sigma < \tau$. (For a
formal definition, see \cite{cms}, pp. 411-412, \cite{som}.)
\end{rem}

\begin{rem}
The definition above is the one introduced in \cite{cms}. For equivalent definitions of fatness, see \cite{ca1},
\cite{ca2}, \cite{pe}, \cite{tu}.
\end{rem}

The idea of the proof of Theorem 1.6. is, in a nutshell, as follows: first build two fat triangulations:
$\mathcal{T}_{1}$ of a certain closed neighbourhood $N_e^\ast$ of the singular set of $\mathbb{H}^{n}/\,G$\,; and
$\mathcal{T}_{2}$ of $(\mathbb{H}^{n}/\,G)\: \backslash \,N_e^\ast$; and then "mash" the two triangulations into a
new triangulation $\mathcal{T}$, while retaining their fatness. The lift $\widetilde{\mathcal{T}}$ of
$\mathcal{T}$ to $\mathbb{H}^{n}$ represents the required $G$-invariant fat triangulation.

Since $(\mathbb{H}^{n}/\,G)\: \backslash \,N_e^\ast$ is an analytical manifold, the existence of $\mathcal{T}_{2}$
is assured by Peltonen's result. In \cite{s2} we showed how to build $\mathcal{T}_{1}$ using a generalization of a
theorem of Munkres (\cite{mun}, 10.6) on extending the triangulation of the boundary of a manifold (with boundary)
to the whole manifold. Munkres' technique also provided us with the basic method of mashing the triangulations
$\mathcal{T}_{1}$ and $\mathcal{T}_{2}$. In this paper we present a more direct, geometric method of triangulating
$N_e^\ast$ and mashing the two triangulations. We already employed this simpler method in \cite{s1}, where we
proved Theorem 1.6. in the case $n = 3$. The original technique used in \cite{s1} for fattening the intersection
of $\mathcal{T}_{1}$ and $\mathcal{T}_{2}$  is, however, restricted to dimension $3$. Therefore here we make
appeal to the method employed in \cite{s2}, which is essentially the one developed in \cite{cms}.

This paper is organized as follows: in Section 2 we present some background on elliptic transformations and  we
show how to choose and triangulate the closed neighbourhood $N_e^\ast$ of the singular set of
$\mathbb{H}^{n}/\,G$, and how to select the "intermediate zone" where the two different triangulations overlap.
Section 3 is dedicated to the main task of fattening the common triangulation. In Section 4 we show how to apply
the main result in the construction of a quasimeromorphic mappings from $\mathbb{H}^n$ to
$\widehat{\mathbb{R}^n}$.

\section{Constructing and Intersecting Triangulations}

\subsection{Elliptic Transformations} Let us first recall the basic definitions and notations: A transformation $f\in Isom(H^n)$, $f\neq Id$ is called {\it elliptic} if $(\exists)\,m\geq 2 $ s.t. $f^{m} = Id$,
and the smallest $m$ satisfying this condition  is called the {\it order} of $f$. In the $3$-dimensional case the
{\it fixed point set} of f, i.e. $ Fix(f) = \{ x\in H^3 | f(x) = x \}$, is a hyperbolic line and will be denoted
by $A(f)$ -- the {\it axis of f}. In dimension $n \geq 4$ the fixed set of an elliptic transformation is a
$k$-dimensional hyperbolic plane, $0 \leq k \leq n-2$\,. The situation is complicated further by the fact that
different elliptics may well have fixed loci of different dimensions (inclusive the degenerate case $n = 0$). If
$A$ is an axes of an elliptic transformation of order $m$ acting upon $\mathbb{H}^3$, then $A$ is called an
$m-axis$. By extension we shall call the $k$-dimensional fixed point set of an elliptic transformation acting upon
$H^{n},\, n \geq 4$\,, an $(m,\,k)$-axis, $0 \leq k \leq n-2$, or just an axis.

\begin{rem}
Since our main interest lies in Kleinian groups acting upon $\mathbb{H}^3$, whose elliptic elements have orders
bounded from above, it should be mentioned that, while for any finitely generated Kleinian group acting on $H^3$
the number of conjugacy classes of elliptic elements is finite (see \cite{fm}), this is not true for Kleinian
groups acting upon $H^n, \,n\geq 4$\,; (for counterexamples, see \cite{fm}, \cite{po} and \cite{h}).
\end{rem}

If the discrete group $G$ is acting upon $\mathbb{H}^n$, then by the discreteness of $G$, there exists no
accumulation point of the elliptic axes in $\mathbb{H}^n$. Moreover, if $G$ contains no elliptics with
intersecting axes, then the distances between the axes are, in general, bounded from bellow. In the classical case
$n = 3$ this bounds are obtained by applying  J\o rgensen's inequality and its corollaries. While no actual
computations in higher dimension of the said bounds are known to us, their feasibility follows from the existence
of various generalizations of J\o rgensen's inequality (see \cite{fh}, \cite{m}, \cite{wa}\,).\footnote{\, For
bounds in dimension $3$ see \cite{bm}, \cite{gm1}, \cite{gm2}.} We defer such computations for further study. The
methods above do not apply for pairs of order two elliptics (see \cite{ah}, \cite{s1}). In this case a
modification of the basic construction will be required (see Section 2.3. below). In the presence of {\it nodes}
i.e. intersections of elliptic axes, the situation is more complicated. Indeed, even in dimension $3$ the actual
computation of the distances between node points has been achieved only relatively recently (see \cite{gm1},
\cite{gm2}, \cite{gmmr})\footnote{\, See also \cite{dm}, \cite{med} and \cite{s3} for a different approach.}.


\subsection{Geometric Neigbourhoods}
To produce the desired closed neighbourhood $N^\ast_e$ of the singular locus of $\mathbb{H}^n/\,G$ and its
triangulation $\mathcal{T}_1$, we start by constructing a standard neighbourhood $N_f = N\big(A(f)\big)$ of the
axes of each elliptic element of $G$ such that $N_f \simeq A(f) \times I^{n-k}$, where $A(f) = \mathbb{S}^k$ and
where $I^{n-k}$ denotes the unit $(n-k)$-dimensional interval. The construction of $N_f$ proceeds as follows: By
\cite{cox}, Theorem $11\cdot23.$ the fundamental region for the action of the stabilizer group of the axes of $f$,
$G_f = G_{A(f)} = \{g \in G\,|\, g(x) = x\}$  is a simplex or a product a simplices. Let $\mathcal{S}_f$ be this
fundamental region. Then we can define the generalized prism\footnote{\, or { \it simplotope} cf. \cite{som}, pp.
113-115.} $\mathcal{S}_f^\perp$, defined by translating $\mathcal{S}_f$ in a direction perpendicular to
$\mathcal{S}_f$, where the translation length is $dist_{hyp}\big(\mathcal{S},A(f)\big)$. It naturally decomposes
into simplices (see \cite{som}, p. 115, \cite{mun}, Lemma 9.4). We have thus constructed an $f$-invariant
triangulation of a prismatic neigbourhood $N_f$ of $A(f)$. We can reduce the mesh of this triangulation as much as
required, while controlling its fatness by dividing it into a finite number of radial strata of equal width
$\varrho = \delta/\kappa_0$, and further partition it into "slabs" of equal hight $h$, where $\delta =
\min\{dist_{hyp}\big(A(f),A(g)\big)\,|\, g \;{\rm elliptic}, \;g \neq f\}$. Henceforth we shall call the
neighberhood thus produced, together with its fat triangulation, a {\it geometric neighbourhood}.

\subsection{Mashing Triangulations}
Since $G$ is a discrete group, $G$ is countable so we can write $G = \{g_j\}_{j \geq 1}$ and let $\{f_i\}_{i \geq
1} \subset G$ denote the set of elliptic elements of $G$. The steps in building the desired fat triangulation are
as follows:

\begin{enumerate}
\item Consider the geometric neighbourhoods \[N_i = N_{i/4} = \{x \in \mathbb{H}^3 \,|\, dist_{hyp}(A_i,x) < \delta/4\}\] and \[N'_i =  N_{i,3/16} = \{x \in \mathbb{H}^n \,|\, dist_{hyp}(A_i,x) <
3\delta/16\}\,,\] with their natural fat triangulations,  where $A_i = A(f_i)$ and where the choice of
"$\delta/4$" instead of "$\delta$" in the definition of the geometric neighbourhoods $N_i$ is dictated by the
following Lemma:

\begin{lem}(\cite{rat})
      Let $X$ be a metric space, and let $\Gamma < Isom(X)$ be a discontinuous group.
      \\Then, for any $x \in X$ and any $r \in (0,\delta /4)$:
      \[ \pi : B(x,r)/\Gamma_{x} \simeq B(\pi(x),\delta/4); \]
      where: $\Gamma_{x}$ is the stabilizer of $x$, $\pi$ denotes the natural projection,
\\ $\delta = d(x,\Gamma(x) \backslash \{x\})$, and where the metric on $X/\Gamma$ is given by
      \[ d_{\Gamma}([\pi(x)],[\pi(y)]) = d(\Gamma(x),\Gamma(y))\,; \,\forall x,y\in  X/\,\Gamma \,.\]
\end{lem}

We also put: $N_e = \bigcup \raisebox{-0.7em}{\hspace{-0.5cm}\tiny $i\in \mathbb{N}$}\!N_{i}$\,.

\item Denote by $\mathcal{T}_e$ the geometric triangulation of $N_e$ described in Section 2.2. above.

\item Consider the following quotients: $N_e^{\ast} =  (\overline{N}_e \cap \mathbb{H}^n)/G$ and $M_c = (\mathbb{H}^n/\,G)\,\backslash \,N_e^{\ast}\\ = (\mathbb{H}^n\, \backslash\, \overline{N}_e)/\,G$.

\item Denote by $\mathcal{T}_p$ the fat triangulation of $M_c$ assured by Peltonen's
Theorem.

\item Consider also the tubes $T_i = N_{i,1/4} \backslash  N_{i,3/16}$, 
and denote $T = \bigcup \raisebox{-0.7em}{\hspace{-0.5cm}\tiny $i\in \mathbb{N}$}\!T_{i}$\,.

\item $T$ is endowed with two triangulations: a natural triangulation $\mathcal{T}_1$ induced by the geometric triangulations of the tubes $T_{i}$, and $\mathcal{T}_2$ that
is inherited from $\mathcal{T}_p$\,.

\item Mashing $\mathcal{T}_1$ and $\mathcal{T}_2$ and ensuring (by a further eventual barycentric type subdivision)
that each $(n-2)$-face of every $n$-simplex is incident to an even number of $n$-simplices, produces a
triangulation $\mathcal{T}_0$\,.

\item Denote by $\mathcal{T}$ the fat triangulation obtained by fattening $\mathcal{T}_0$\,.

\item Let $\widetilde{\mathcal{T}}$ be the lift of $\mathcal{T}$ to $\mathbb{H}^n$. Then $\widetilde{\mathcal{T}}^\sharp = \widetilde{\mathcal{T}} \cup \mathcal{T}_e$  represents the desired $G$-invariant
fat triangulation of $\mathbb{H}^n$.

\end{enumerate}

\begin{rem}
The construction above applies only in the case of non-intersecting elliptic axes. In the case when there exist
intersecting elliptic axes, the following modification of our construction is required: instead of $\delta$ one
has to consider $\delta^* = \min{(\delta, \delta_0)}$, where $\delta_0$ represents the minimal distance between
nodes.
\end{rem}

\begin{rem}
In the case when all the elliptic transformations are half-turns, since, as we have seen, no minimal distance
between the axes can be computed in this case. However, by the discreetness of $G$ it follows that there is no
accumulation point of the axes in $\mathbb{H}^3$. Let $D = \{d_{ij}\,|\, d_{ij} = dist_{hyp}(A_i,A_j)\}$ denote
the set of mutual distances between the axes of the elliptic elements of $G$. Then, since $G$ is countable, so
will be $D$, thus $D = \{d_k\}_{k \geq 1}$. Then  the set of neighbourhoods $N_e^{\natural} = \bigcup
\raisebox{-0.7em}{\hspace{-0.5cm}\tiny $k \in \mathbb{N}$}\!N_{k}^{\natural} = \bigcup
\raisebox{-0.7em}{\hspace{-0.5cm}\tiny $k \in \mathbb{N}$}\!\{x \in \mathbb{H}^n \,|\, dist_{hyp}(A_k,x) <
\delta/4k\}$ will constitute a proper geometric neigbourhood of $A_G =
\bigcup\raisebox{-0.7em}{\hspace{-0.5cm}\tiny $i \in \mathbb{N}$}\!A_i$. The fatness of the simplices of the
geometric triangulation of $A_G$ can be controlled, as before, by a proper choice of $h$ and $\varrho$.
\end{rem}

\section{Fattening Triangulations}

First let us establish some definitions and notations:
\\ Let $K$ denote a simplicial complex, let $K' < K$ denote a subcomplex of $K$ and let $\sigma \in K$ denote the simplices of $K$.

\begin{defn}
Let $\sigma_i \in K$, $dim\,\sigma_i = k_i$, $i=1,2$; s.t. $diam\,\sigma_1 \leq diam\,\sigma_2$. We say that
$\sigma_1, \sigma_2$ are $\delta${\em-transverse} iff
\\ (i) $dim(\sigma_1 \cap \sigma_2) = \max(0,k_1+k_2-n)$;
\\ (ii) $0 < \delta < \measuredangle(\sigma_1,\sigma_2)$;
\\ and if $\sigma_3 \subset \sigma_1,\,\sigma_4 \subset \sigma_2$, s.t. $dim\,\sigma_3 + dim\,\sigma_4 < n =
dim\,K$, then
\\ (iii) $dist(\sigma_3, \sigma_4) > \delta\cdot\eta_1$.
\\ In this case we write: $\sigma_1 \pitchfork_{\delta} \sigma_2$.
\end{defn}

We begin by triangulating and fattening the intersection of two individual simplices belonging to the two given
triangulations, respectively. Given two closed simplices $\bar{\sigma_1}, \bar\sigma_2$, their intersection (if
not empty) is a closed, convex polyhedral cell: $\bar\gamma = \bar\sigma_1 \cap \bar\sigma_2$. One canonically
triangulates  $\bar\gamma$ by using the {\em barycentric subdivison} $\bar\gamma^\ast$ of $\bar\gamma$, defined
inductively upon the dimension of the cells of $\partial\gamma$  in the following manner: for each cell $\beta
\subset \partial\gamma$, choose an interior point $p_{\beta} \in int\,\beta$ and construct the join
$J(p_\beta,\beta_i),\; \forall \beta_i \subset \partial \beta$.\footnote{\, If $dim\,\beta = 0$ or $dim\,\beta =
1$, then $\beta$ is already a simplex.}

We first show that if the simplices are fat and if they intersect $\delta$-transversally, then one can choose the
points s.t. the barycentric subdivision $\bar\gamma^\ast$ will be composed of fat simplices. More precisely, we
prove the following Proposition:

\begin{prop}
Let $\sigma_1, \sigma_2 \subset \mathbb{R}^m$, where $m = \max(dim\,\sigma_1,dim\,\sigma_2)$, s.t.
\\ $d_1 = diam\,\sigma_1 \leq d_2 = diam\,\sigma_2$, and s.t. $\sigma_1, \sigma_2$ have common fatness $\varphi_0$.
\\ If $\sigma_1 \pitchfork_{\delta} \sigma_2$, then there exists $c = c(m,\varphi_0,\delta)$

\begin{enumerate}
\item If $\sigma_3 \subset \bar\sigma_1,\,\sigma_4 \subset \bar\sigma_2$ and if $\sigma_1 \cap \sigma_2 \neq \emptyset$, then $\sigma_1 \cap \sigma_2 =
\gamma_0$ is an ($k_3+k_4-m$)-cell, where $k_3 = dim\,\sigma_3, k_4 = dim\,\sigma_4$ and:

\begin{equation}
Vol_{eucl}(\gamma_0) \geq c\cdot d_1^{k_3+k_4-m}\;.
\end{equation}

\item $\forall\, \gamma_0$ as above, $\exists\, p \in \gamma_0$, s.t.
\begin{equation}
dist(p,\partial\gamma_0) > c\cdot d_1\;.
\end{equation}

\item If the points employed in the construction of $\bar\gamma^\ast$ satisfy the condition {\rm (3.2)} above, then each
$l$-dimensional simplex $\tau \in \gamma^\ast$ satisfies the following
\\inequalities:

\begin{equation}
\varphi_l \geq Vol_{eucl}(\tau)/d_1^{\,l} \geq  c\cdot d_1\;.
\end{equation}

\end{enumerate}

\end{prop}

\begin{proof}

First, consider the following remarks:

\begin{rem}
 The following sets are compact:
\\ $S_1 = \{\sigma_1\,|\, diam\,\sigma_1 = 1\,, \varphi(\sigma_1) \geq \varphi_0\}$, $S_2 = \{\sigma_2\,|\, diam\,\sigma_2 = 2(1+\delta)\,, \varphi(\sigma_2) \geq
\varphi_0\}$,

$S(\phi_0,\delta) \subset S_1 \cap S_2\,, \; S(\phi_0,\delta) = \{(\sigma_1, \sigma_2)\,|\, \exists v_0, {\rm
s.t.} v_0 \in \sigma_1, \forall \sigma_1 \in S_1 \cap S_2\}$.
\end{rem}

\begin{rem}
There exists a constant $c(\varphi)$ s.t. $\mathcal{S} = \mathcal{S}'$, where
\\ $\mathcal{S} = \{\sigma_1 \cap \sigma_2\,|\, diam\,\sigma_2 \leq d_2\},\;  \mathcal{S}' = \{\sigma_1 \cap \sigma_2\,|\, diam\,
c(\varphi)(1+\delta)d_1\}$,
\\ i.e. the sets of all possible intersections remains unchanged under controlled dilations of one of the families
of simplices.
\end{rem}

Now, from the fact that $\sigma_1 \pitchfork_{\delta} \sigma_2$ it follows that $\sigma_3 \cap \sigma_4 \neq
\emptyset \Leftrightarrow \bar{\sigma}_3 \cap \bar{\sigma}_4 \neq \emptyset$ (see \cite{cms}, p. 436). Therefore,
the function $Vol_{eucl}(\gamma_0)$ attains a positive minimum, as a positive, continuous function defined on the
compact set $\bar{\sigma}_3 \cap \bar{\sigma}_4$, thus proving the first assertion of the proposition.

Let be $\gamma$ be a $q$-dimensional cell, and let $\beta$ be a face of $\partial\gamma$. Then:

\begin{equation}
Vol_{eucl}(\beta) \leq d_1^{p}\;.
\end{equation}

Choose $p \in \gamma$, such that $\rho = dist(p,\partial\gamma) = \max\{dist(r,\partial\gamma)\,|\,r \in
\gamma\}$. Then, if $\beta = \beta^j$ denotes a $j$-dimensional face of $\partial\gamma$, we have that:

\begin{equation}
\gamma \subset \bigcup_{\beta^j \subset \partial\gamma}\!\!N_{\rho}(\beta^j);
\end{equation}

where: $N_{\rho}(\beta^j) = \{r\,|\, dist(r,\beta^j) \subseteq \rho\}$. But:

\begin{equation}
Vol_{eucl}(\beta^j \cap \gamma) \leq c\cdot\rho^{q-j}\cdot Vol_{eucl}(\beta^j),
\end{equation}

for some $c' = c'(q)$.
\\ Moreover, the number of faces $\sigma_3 \cap \sigma_4$ of $\gamma$ is at most
$2^{dim\,\sigma_1+dim\,\sigma_2+2}$, where $\sigma_1, \sigma_2$ are as in Remark 3.4. and  $dim\,\sigma_3 \leq
dim\,\sigma_1,\;  dim\,\sigma_4 \leq dim\,\sigma_2$.
\\Thus (3.4) in conjunction with (3.6) imply that there exists $c_1 = c_1(m,\varphi_0,\delta)$, such that:

\begin{equation}
c_1d_1^q \leq \sum_{j=0}^{q-1}\rho^{q-j}d_1^j.\,
\end{equation}

and (3.2) follows from this last inequality.

The last inequality follows from (3.2) and (3.3) by induction.

\end{proof}

Next we show that given two fat Euclidian triangulations that intersect $\delta$-transver-sally, then one can
infinitesimally move any given point of one of the triangulations s.t. the resulting intersection will be
$\delta'$-transversal, where $\delta'$ depends only on  $\delta$, the common fatness of the given triangulations,
and on the displacement length. More precisely one can show that the following results holds:

\begin{prop}
Let $K_1, K_2 \subset \mathbb{R}^n$ be $n$-dimensional simplicial complexes, of common fatness $\varphi_0$ and
$d_1 = diam\,\sigma_1 \leq d_2 = diam\,\sigma_2$. Let $v_0 \in K_1$ be a $0$-dimensional simplex of $K_1$.
Consider the complex $K^*_1$ obtained by replacing $v_0$ by $v^*_0 \in R^n$ and keeping fixed the rest of the
$0$-dimensional vertices of $K_1$ fixed. Consider also $L_2 < K_2,\; L_2 = \{\sigma \in K_2\,|\, \sigma \cap
B(v_0,2d_1) \neq \emptyset\}$.
\\ Then,  if there exists $k$ s.t. all the $k$-simplices $\tau \subset \partial St(v_0)$ are $\delta$-transversal
to $L_2$, there exist $\varphi_0, \delta, \varepsilon > 0$, $\delta^\ast =
\delta^\ast(\varphi_0,\delta,\varepsilon)$ and there exists $v^*_0$ s.t. $dist(v_0,v^*_0) < \varepsilon \cdot d_1$
s.t.

\begin{equation}
\tau^* \pitchfork_{\delta^{\ast}} L_2\,; \;\forall\, \tau \subset St(v^*_0)\, \backslash\, \partial St(v^*_0),\;
dim\, \tau^* = k+1.
\end{equation}

\end{prop}

\begin{proof}
Let $N(r) = |\{ \sigma \in K_1\,|\, \sigma \subset B_r(v_0)\}|$. Then there exists a constant $c_n$ s.t. $N(r)
\leq \frac{c_n}{\varphi_0}(\frac{\varepsilon}{d_1})^n$. It follows that the set $St(v_0)$ is compact, since there
are at most $\frac{c_n}{\varphi_0}$ possible edge lengths, which can take values in the interval
$[d_1\varphi_0,d_1]$.\footnote{\,i.e. the number of possible combinatorial structures on $St(v_0)$ depends only on
$\varphi_0$.} Therefore if a $D^\ast$ satisfying (3.8) exist, it depends only on $\varphi_0,\; \delta$ and
$\varepsilon$ (and not on $K_1, K_2$).
\\ Let $\sigma_1, \ldots , \sigma_{l_1}$ and $\tau_1, \ldots, \tau_{l_2}$ be orderings of the simplices of $L_2$ and of the $k$-simplices of $\partial St(v_0)$, respectively.
Then, by \cite{cms}, Lemma 7.4, there exists $\varepsilon_{1,1}$ and exists $v_{1,1}, \; d(v_{1,1},v_0) =
\varepsilon_{1,1}$, such that the hyperplane $\Pi(v_{1,1},\tau_1)$ determined by $v_1$ and by $\tau_1$ is
transversal to $\sigma_1$. By replacing $\tau_1$ by $\tau_2$ and $v_0$ by $v_{1,1}$ we obtain $v_{1,2}$ and
$\varepsilon_{1,2}$ s.t. $\Pi(v_{1,2},\tau_1) \pitchfork \sigma_2$. (See Fig. 1.)

\begin{figure}[h]
\begin{center}
\includegraphics[scale=0.4]{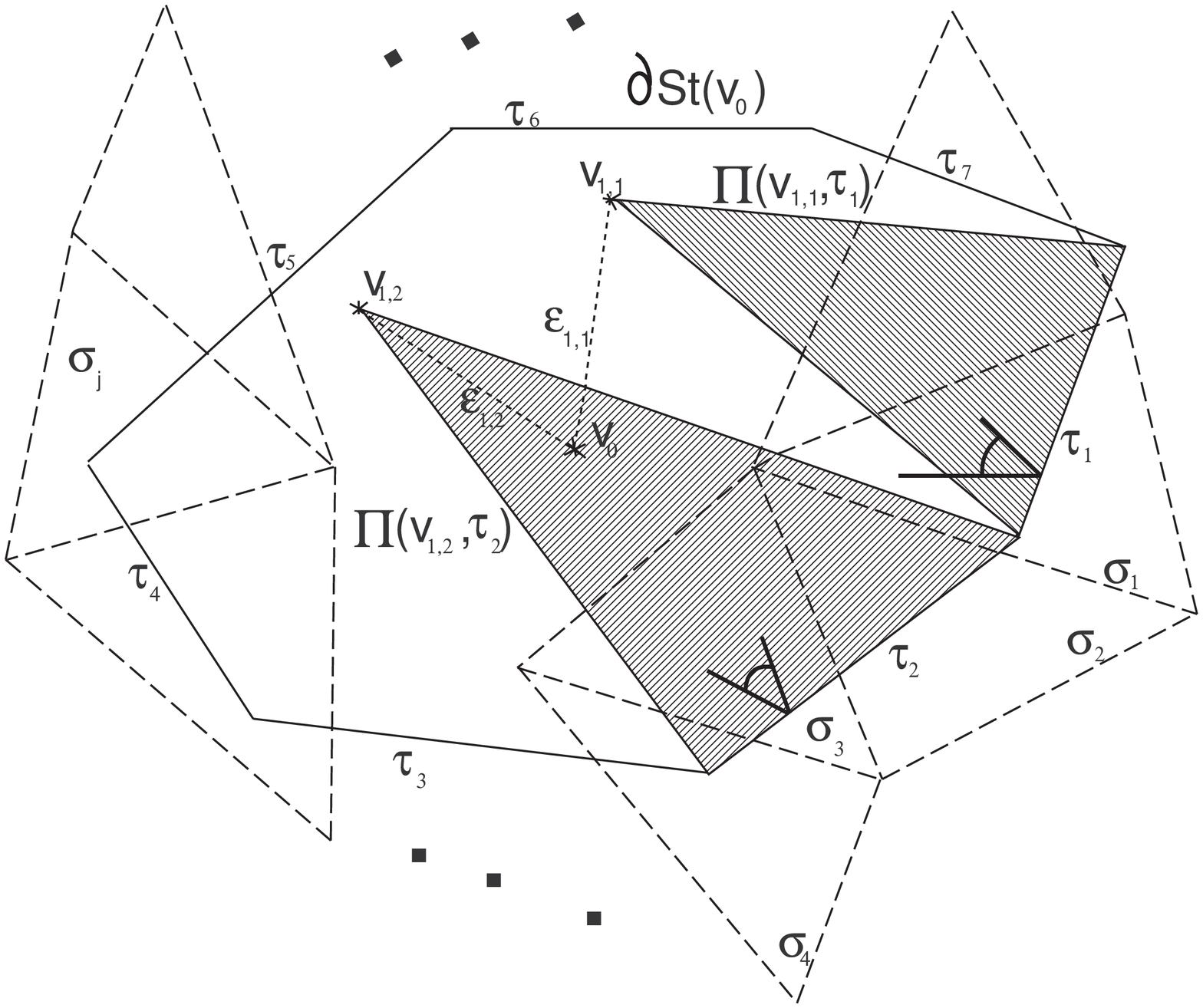}
\end{center}
\caption{}
\end{figure}

Moreover, by choosing $\varepsilon_{1,2}$ sufficiently small, one can ensure that $\Pi(v_{1,2},\tau_1) \pitchfork
\sigma_1$, also. Repeating the process for $\tau_3, ..., \tau_{l_2}$, one determines a point  $v_{1,l_2}$ such
that $\Pi(v_{1,l_2},\tau_j) \pitchfork L_2\; j = 3,\ldots,l_2$. In the same manner and by choosing at each stage
an $\varepsilon_{i,j}$ small enough, one finds points $v_{i,j}$ s.t. $\Pi(v_{i,j},\tau_i) \pitchfork L_2,\; i = 1,
\ldots ,l_1\,,\; j = 1,\dots,l_2$. Then $v^*_0 = v_{l_1,l_2}$ satisfies: $\Pi(v_0^\ast,\tau_{j}) \pitchfork L_2$,
$j = 1,\ldots,,l_2$.

\end{proof}

We are now prepared to prove the main result of this section namely:

\begin{thm}
Let $M^n$ be an orientable $n$-manifold and let $\mathcal{T}_1,\mathcal{T}_2$ be two fat triangulations of open
sets $U_1,U_2 \subset M^n,\,$ $U_1 \cap U_2 \neq \emptyset$, having common fatness $\geq \varphi_0$, and such that
$\mathcal{T}_1 \cap \mathcal{T}_2 \neq \emptyset$. Then there exist fat triangulations
$\widetilde{\mathcal{T}'_1},\mathcal{T}'_2$ and there exist open sets $U \subset  U_1 \cap U_2 \subset V$, such
that
\begin{enumerate}
\item $({\mathcal{T}'_1} \cap \mathcal{T}'_2) \cap (U_i \setminus V) = \mathcal{T}_i\,,\; i=1,2\,;$
\item $({\mathcal{T}'_1} \cap \mathcal{T}'_2) \cap U = \mathcal{T};$
\\\hspace*{-1.3cm} where
\item $\mathcal{T}$ is a fat triangulation of $U$.
\end{enumerate}
\end{thm}

\begin{proof}
Let $K_1,\,K_2$ denote the underlying complexes of $\mathcal{T}_1,\mathcal{T}_2$, respectively. By considerations
similar to those of Proposition 3.5. it follows that given $\varphi_0
> 0$, there exists $d(\varphi_0) > 0$ such that given a $k$-dimensional simplex $\sigma \subset \mathbb{R}^n$,
$diam(\sigma) = d_1$ has fatness $\varphi_0$, than translating each vertex of $\sigma$ by a distance
$d(\varphi_0)\cdot d_1$ renders a simplex of fatness $\geq \varphi_0/2$\,. Also, it follows that given
$\varphi_0,\delta
> 0$, exists $\delta(\varphi_0,\delta)$ satisfying the following condition: if every vertex $u \in \sigma \subset
K$ is replaced by a vertex $u'$ s.t. $dist(u,u') \leq \delta(\varphi_0,\delta) \cdot d_1$, then the resulting
simplex $\sigma'$ is $\pitchfork_{\delta/2}$-transversal to $K$; for any $n$-dimensional simplicial complex $K$ of
fatness $\varphi_0$ and such that $diam\,\sigma = d_2 \geq d_1$ .
\\ Let $v_0 \in U_1 \cap U_2$. Define the following subcomplexes of $K_1,\,K_2$, respectively:
\[L_2 = \{\bar{\sigma} \subset K_2 \,|\, \bar{\sigma} \subset B_\varepsilon(v_0),\, d_1 \leq dist(\bar{\sigma},\partial B_\varepsilon(v_0)) \leq d_2\}\]
\[M_2 = \{\sigma \subset K_2 \,|\, \bar{\sigma} \subset \bar{\tau} \subset B_\varepsilon(v_0),\, dim\,\tau = n,\,\bar{\sigma} \cap L_2 \neq \emptyset\}\]
\[\hspace*{-3.6cm}L_1 = \{\bar{\sigma} \subset K_1 \,|\, dist(\bar{\sigma},L_2) \leq d_2\}\]
\[M_1 = \{\sigma \subset K_1 \,|\, \bar{\sigma} \subset \tau \subset B_\varepsilon(v_0),\, dim\,\tau = n,\, \tau \cap L_1 \neq \emptyset\}\]
(See Fig. 2.)

\begin{figure}[h]
\begin{center}
\includegraphics[scale=0.35]{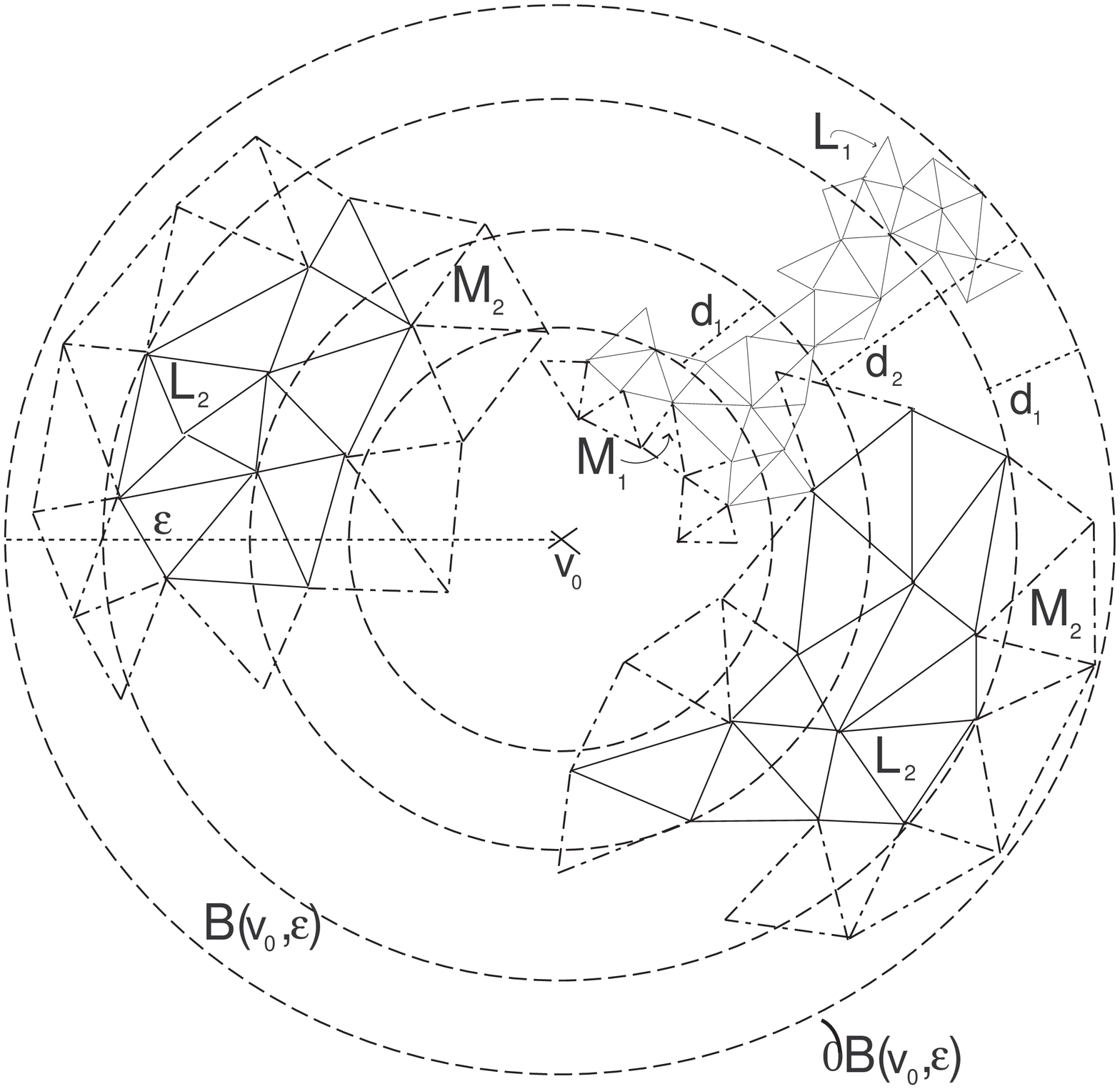}
\end{center}
\caption{}
\end{figure}

Consider an ordering $v_1, \ldots, v_p$ of the vertices of $L_1$. It follows from Proposition \nolinebreak[4]3.5.
that, if all the vertices of $L_1$ are moved by at most $t_0$, where

\begin{equation}
t_0 = \frac{d_1}{n}\min{\big\{\frac{1}{2},d(\varphi_0)\big\}}\,,
\end{equation}

then there exists

\begin{equation}
\delta^\ast_0 = \delta_0^\ast(\varphi_0,1,\frac{t_0}{d_1})\,,
\end{equation}

such that

\begin{equation}
\mathcal{S}^0(L_{1,0}) \pitchfork_{\delta_{0}^\ast} K_2\,,
\end{equation}

where $\mathcal{S}^0(L_{1,0})$ denotes the $0$-skeleton of $L_{1,0}$.

Now define inductively

\begin{equation}
t_i =\frac{d_1}{n}\min{\bigg\{\frac{1}{2},d(\varphi_0),\delta\Big(\frac{\varphi_0}{2},\frac{\delta^\ast_0}{2}\Big)
,\ldots,\delta\Big(\frac{\varphi_0}{2},\frac{\delta^\ast_{i-1}}{2}\Big) \bigg\}}\,,
\end{equation}

where

\begin{equation}
\delta^\ast_{i} = \delta^\ast_{i}\Big(\varphi_0,\delta^\ast_{i-1},\frac{t_i}{d_1}\Big)\,;\; i=1,\ldots,n-1.
\end{equation}

Then $t_0 \geq t_1 \geq \ldots \geq t_{n-1}$.
\\ Moving  each and every vertex of $L_1$ by a distance $\leq t_i,\;i=1,\ldots,n$, renders complexes
$L_{1,1},\ldots,L_{1,n-1}$ s.t.
\begin{enumerate}
\item $L_{1,i}\cap\big(B_\varepsilon(v_0) \setminus \, M_1 \big) \equiv L_1$,
\item $L_{1,i}$ are $\varphi_0$-fat; $i=0,\ldots,n-1$.
\end{enumerate}

By inductively applying Proposition 3.5. it follows that

\begin{equation}
\mathcal{S}^i(L_{1,i}) \pitchfork_{\delta_i^\ast} K_2\,,
\end{equation}

Where $\mathcal{S}^i(L_{1,i})$ denotes the $i$-skeleton of $L_{1,i}$\,.
\\Moreover,

\begin{equation}
\mathcal{S}^i(L_{1,j}) \pitchfork_{\delta_i^\ast/2} K_2\,, \forall j > i.
\end{equation}

It follows that

\begin{equation}
L_{1,n-1} \pitchfork_{\delta^\bigstar} K_2\,,
\end{equation}

where

\begin{equation}
\delta^\bigstar = \frac{1}{2}\min\{\delta_0^\ast,\ldots,\delta_i^\ast\}\,.
\end{equation}

By Proposition 3.2. the barycentric subdivision of $L_{1,n-1} \cap L_2$ is fat. We extend it to a fat subdivision
of $M_2$ in the following manner: given a simplex $\sigma \subset  M_2 \setminus L_2$, subdivide $\sigma$ by
constructing all the simplices with vertices $v_i$, where $v_i$ is either the vertex of a simplex $\sigma \subset
M_2 \setminus L_2,\; \bar{\sigma_i} \cap L_2 \neq \emptyset$, or it is a vertex of a closed simplex $\bar{\sigma}$
of the barycentric subdivision of $L_{1,n-1} \cap L_2$\,, such that $\bar{\sigma} \subset \partial L_2 \cap
M_2\,,\; i = 1,\ldots,k_0$. The triangulation $\widetilde{K}_2$ thus obtained is a fat extension of $\overline{K_2
\setminus\,M_2 }$. 

In an analogous  manner one constructs a similar fat extension $\widetilde{K}_1$ of $\overline{K_1
\setminus\,M_2}$.

\end{proof}

Now let $M^n= M_c$ and let $\mathcal{T}_1$, $\mathcal{T}_2$ be the fat chessboard triangulations constructed in
2.3.(6)--2.3.(7)\,. Then the local fat triangulation obtained in Theorem 3.6. above extends globally to the fat
triangulation of $\mathcal{T}_1 \cap \mathcal{T}_2$, by applying Lemma 10.2 and Theorem 10.4 of \cite{mun}. This
provides us with the sought fat chessboard triangulation $\mathcal{T}$ of $M_c$. Thus
$\widetilde{\mathcal{T}}^\sharp = \widetilde{\mathcal{T}} \cup \mathcal{T}_e$ represents the required
$G$-invariant fat chessboard triangulation of $\mathbb{H}^n$.

\section{The Existence of Quasimeromorphic Mappings}

The technical ingredient in Alexander's trick is the following Lemma:
\begin{lem} (\cite{ms1}, \cite{pe})
Let $M^n$ be an orientable $n$-manifold, let $\mathcal{T}$ be a chessboard fat triangulation of $M^n$ and let
$\tau,\sigma \in \nolinebreak[4]\mathcal{T},\; \tau = (p_1,\dots,p_n), \,\sigma = (q_1,\dots,q_n)$; and denote
$|\tau| = \tau \cup int\,\tau$.
\\ Then there exists a orientation-preserving homeomorphism $h = h_{\tau}: |\tau| \rightarrow \widehat{\mathbb{R}^{n}}$
s.t.
\begin{enumerate}
\item $h(|\tau|) = |\sigma|$, \,if\; $det(p_1,\dots,p_n) > 0$
\\ and
\\ $h(|\tau|) = \widehat{\mathbb{R}^{n}} \setminus| \sigma|$, \,if\; $det(p_1,\dots,p_n) < 0$.
\item $h(p_i) = q_i, \; i=1,\ldots,n.$
\item $h|_{\partial|\sigma|}$ is a $PL$ homeomorphism.
\item $h|_{int|\sigma|}$ is quasiconformal.
\end{enumerate}
\end{lem}

\begin{prf}
Let $\tau_0 = (p_{0,1},\dots,p_{0,n})$ denote the equilateral $n$-simplex inscribed in the unit sphere
$\mathbb{S}^{n-1}$. The {\it radial linear stretching} $\varphi: \tau \rightarrow \overline{B^n}$ is onto and
bi-lipschitz (see \cite{ms2}). Moreover, by a result of Gehring and V\"{a}isal\"{a}, $\varphi$ is also
quasicomformal (see \cite{v}). We can extend $\varphi$ to $\widehat{\mathbb{R}^n}$ by defining $\varphi(\infty) =
\infty$. Let $J$ denote the reflection in the the unit sphere $\mathbb{S}^{n-1}$ and let $h_0: |\sigma|
\rightarrow |\tau|$ denote the orientation-reversing $PL$ mapping defined by: $h(p_i) = q_i, \; i=1,\ldots,n$.
Then $h = \varphi^{-1} \circ J \circ \varphi \circ h_0$ is the required homeomorphism.
\end{prf}

The Existence Theorem  of quasimeromorphic mappings now follows immediately:
\vspace*{0.2cm}\\ \begin{prf} {\bf of Theorem 1.6.} Since the orders of the elliptic transformations are uniformly
bounded, so will be the fatness of the simplices of $\mathcal{T}_1$ -- the geometric triangulation  of $T$. Let
$\widetilde{\mathcal{T}}^\sharp$ be the $G$-invariant fat chessboard triangulation of $\mathbb{H}^n$ constructed
above. Let $f:\mathbb{H}^n\rightarrow \widehat{\mathbb{R}^n}$ be defined by: $f|_{|\sigma|} = h_{\sigma}$, where
$h$ is the homeomorphism constructed in the Lemma above. Then $f$ is a local homeomorphism on the $(n-1)$-skeleton
of $\widetilde{\mathcal{T}}^\sharp$ too, while its branching set $B_f$ is the $(n-2)$-skeleton of
$\widetilde{\mathcal{T}}^\sharp$. By its construction $f$ is quasiregular. Moreover, given the uniform fatness of
the simplices of triangulation $\widehat{\mathcal{T}}^\sharp$, the dilatation of $f$ depends only on the dimension
$n$.
\end{prf}


\end{document}